\documentclass[preprint,12pt]{elsarticle}

\usepackage[tbtags]{amsmath}
\usepackage{bm}
\usepackage{amsfonts}
\usepackage{graphics}
\usepackage{graphicx}
\usepackage{psfrag}
\usepackage{eurosym}
\usepackage{pstricks}
\usepackage{array}
\usepackage{shortvrb}
\usepackage{epsf}
\usepackage{graphicx}
\usepackage{rotating}
\usepackage{float}
\usepackage{color}
\usepackage{multirow, makecell}
\usepackage{float}
\usepackage{colortbl}
\usepackage{ifpdf}
\usepackage{multicol}

\usepackage{hyperref}
\hypersetup{colorlinks=true, linkcolor=blue, anchorcolor=red, citecolor=blue, filecolor=red, urlcolor=red, pdfauthor=author}

\usepackage{geometry}
\usepackage{alphalph}
\usepackage{etoolbox}
\usepackage{pdfpages}

\usepackage[T1]{fontenc}

\usepackage{caption}
\usepackage{subcaption}
\newcommand{\comment}[1]{}

\patchcmd{\subequations}{\alph{equation}}{\alphalph{\value{equation}}}{}{}
\def \Noin {{\hskip 3pt \rm /\kern -9pt \in\hskip 1pt}}

\def \R {{\rm I\kern -2.2pt R\hskip 1pt}}

\allowdisplaybreaks

\begin{document}
%\includepdf[pages=-, fitpaper]{Portada_WP22_06.pdf}

\begin{frontmatter}

\title{On data-driven chance constraint learning for mixed-integer optimization problems}

\author[inst1]{Antonio Alc\'antara\corref{cor1}}
\ead{antalcan@est-econ.uc3m.es}
\author[inst1,inst2]{Carlos Ruiz}
\ead{caruizm@est-econ.uc3m.es}
\cortext[cor1]{Corresponding author}

\affiliation[inst1]{organization={Department of Statistics, University Carlos III of Madrid},%Department and Organization
            addressline={Avda. de la Universidad, 30}, 
            city={Legan\'es},
            postcode={28911}, 
            state={Madrid},
            country={Spain}}

\affiliation[inst2]{organization={UC3M-BS Institute for Financial Big Data (IFiBiD), University Carlos III of Madrid},%Department and Organization
            addressline={Avda. de la Universidad, 30}, 
            city={Legan\'es},
            postcode={28911}, 
            state={Madrid},
            country={Spain}}

\begin{abstract}
When dealing with real-world optimization problems, decision-makers usually face high levels of uncertainty associated with partial information, unknown parameters, or complex relationships between these and the problem decision variables. In this work, we develop a novel Chance Constraint Learning (CCL) methodology with a focus on mixed-integer linear optimization problems which combines ideas from the chance constraint and constraint learning literature. Chance constraints set a probabilistic confidence level for a single or a set of constraints to be fulfilled, whereas the constraint learning methodology aims to model the functional relationship between the problem variables through predictive models. One of the main issues when establishing a learned constraint arises when we need to set further bounds for its response variable: the fulfillment of these is directly related to the accuracy of the predictive model and its probabilistic behaviour. In this sense, CCL makes use of linearizable machine learning models to estimate conditional quantiles of the learned variables, providing a data-driven solution for chance constraints. An open-access software has been developed to be used by practitioners. Furthermore, benefits from CCL have been tested in two real-world case studies, proving how robustness is added to optimal solutions when probabilistic bounds are set for learned constraints.
\end{abstract}


\end{frontmatter}

\section{Introduction}
\label{sec:intro}

Recent digitization and automation within our society have not only brought big amounts of data available for decision-makers, but also high levels of uncertainty inherent to the information. From assessing climate change to dealing with the evolution of a pandemic \cite{ertem2021decision}, decision-makers have to set optimal solutions in an environment full of uncertainty.

In the operations research field, there exist certain standard approaches that are typically employed to solve problems where uncertainty and/or big amounts of external data are present. In particular, one of the most employed ones is stochastic programming, where the optimization model incorporates the estimated probability distribution of the uncertain parameters \cite{birge2011introduction}. In particular, stochastic programming considers the following problem:
\begin{subequations}\label{eq:sto_pro}
\begin{align}
\underset{x \in \mathcal{X}}{\min} \;& \mathbb{E}\left[ c(x;\xi) \right] \\
\text{s.t.}& \quad g(x;\xi)\leq 0 \label{eq:sto_pro_cons1}
\end{align}
\end{subequations}

\noindent where $x \in \mathcal{X} \subset \mathbb{R}^{d_x}$ represents the decision variables, $\xi$ is a random vector valued in $\mathbb{R}^{d_{\xi}}$, $c(x;\xi) : \mathbb{R}^{d_x} \times \mathbb{R}^{d_{\xi}} \rightarrow \mathbb{R}$ is the cost function, $g(x;\xi): \mathbb{R}^{d_x} \times \mathbb{R}^{d_{s}} \rightarrow \mathbb{R}^{d_{g}}$ are the constraints, and $\mathbb{E}\left[\cdot \right]$ represents the expected value over the cost function.

Some recent works have focused on upgraded versions of (\ref{eq:sto_pro}) that leverage the inclusion of external information and predictive models within the problem-solving methodology. From the classical ``Predict and Optimize'' strategy, researchers have moved on to frameworks where prediction and prescriptions intersect. For example, in \cite{elmachtoub2022smart}, authors address an integrated approach to find functions $f$ to model other decision variables $y$ influenced by the uncertainty $\xi$ that also generate good prescriptions. The so-called Predictive to Prescriptive two-step approach is proposed in \cite{bertsimas2020predictive} to manage the trade-off between prediction and prescription performance.

An alternative for dealing with uncertainty is the employment of robust optimization methodologies \cite{ben2009robust}. Some examples of robust optimization rely on renewable-energy integration in microgrids \cite{lekvan2021robust}, dealing with risks and improving the sustainability of supply chains \cite{lotfi2021robust, govindan2018advances}, or establishing robust R\&D budget allocation \cite{jang2019decision}. Historically, one of the most employed approaches in the literature has been Wald's maximin model for non-probabilistic robust optimization problems. A simple model is presented in (\ref{eq:robust_pro}), where $\xi$ is a random vector valued in $\mathbb{R}^{d_{\xi}}$, that this time is contained in the uncertainty set $\mathcal{U}(x)$.
\begin{subequations}\label{eq:robust_pro}
\begin{align}
\underset{x \in \mathcal{X}}{\max} \underset{\xi \in \mathcal{U}(x)}{\min} \;&  c(x;\xi) \\
\text{s.t.}& \quad g(x;\xi)\leq 0 \label{eq:robust_pro_cons1}
\end{align}
\end{subequations}

The adequate modelling of the uncertainty set $\mathcal{U}(x)$ has been a major research problem for the last years. In particular, data-driven uncertainty sets may allow us to obtain probabilistic guarantees associated with the optimal solutions to our robust optimization problems \cite{bertsimas2018data}. Several well-known uncertainty sets have been explored \cite{ben2009robust}. For example, a Markov set contains all distributions perfectly defined by its mean and support. In a Chebyshev set, distributions are presented with bounds on the first and second-order moments, whereas in a Huber set all contained distributions present a known upper bound on the expected Huber loss function.

Finally, some applications need to explicitly model and limit the uncertainty associated with the satisfaction of a single or a set of constraints by using probabilities \cite{guo2021chance}. This is the main approach we will address in this work, the so-called chance constraint approach. We will start by introducing the general formulation of a chance-constrained optimization problem (\ref{eq:model_intro1_cons1}). For a fixed safety probability level $\alpha \in (0,1)$, we can replace constraint (\ref{eq:sto_pro_cons1}) by:
%
%\begin{subequations}\label{eq:model_intro1}
%\begin{align}
%\underset{x \in \mathcal{X}}{\min} \;& \mathbb{E}\left[ %c(x;\xi) \right] \\
%\text{s.t.}& \quad \mathbb{P}\left(g(x;\xi)\leq 0\right) \geq \alpha \label{eq:model_intro1_cons1}
%\end{align}
%\end{subequations}
%
\begin{equation}\label{eq:model_intro1_cons1}
\mathbb{P}\left(g(x;\xi)\leq 0\right) \geq \alpha 
\end{equation}

Chance constraint (\ref{eq:model_intro1_cons1}) imposes that $g(x;\xi)\leq 0$ will be satisfied with a probability not below $\alpha$. Traditionally, one approach for addressing chance constraints has been scenario generation. Therefore, the real (or estimated) probability distribution of the uncertain parameters is discretized, and scenarios are generated from this discretization. Then, the big-M approach is used jointly with binary variables to check and assert that at least $\alpha\%$ of the binary variables are set to one. This approach is also known as Sample Average Approximation (SAA) \cite{kuccukyavuz2021chance}. One of the main disadvantages of this proposal is the computational cost increment due to the addition of binary variables for each scenario. 

However, chance constraint (\ref{eq:model_intro1_cons1}) can also be written equivalently as follows:
\begin{equation}\label{eq:quant}
Q_{\alpha}(g(x;\xi)) \leq 0
\end{equation}

\noindent where $Q_{\alpha}(\cdot)$ represents the quantile function associated with probability $\alpha$. This indicates that $\alpha\times 100\%$ of the values of a random variable will lay under the value of this quantile. $Q_{\alpha}(\cdot)$ is also considered in the literature as the Value-at-Risk (VaR) for a certain random variable at probability level $\alpha$.

We can also find in the literature a related approach that employs the buffered failure probability method \cite{rockafellar2010buffered} to transform (\ref{eq:model_intro1_cons1}) into:

%\begin{subequations}\label{eq:model_intro2}
%\begin{align}
%\underset{x \in \mathcal{X}}{\min} \;& \mathbb{E}\left[ %c(x;\xi) \right] \\
%\text{s.t.}& \quad \overline{Q}_{\alpha}(g(x;\xi)) \leq 0  %\label{eq:model_intro2_cons1}
%\end{align}
%\end{subequations}
%
\begin{equation}\label{eq:superq}
\overline{Q}_{\alpha}(g(x;\xi)) \leq 0 
\end{equation}

\noindent where $\overline{Q}_{\alpha}$ represents the superquantile function associated with probability $\alpha$, also known as the Conditional-Value-at-Risk (CVaR). This superquantile function can be simply described as the expectation of the $(1-\alpha)$-tail distribution for a random variable. Regarding its relationship with the quantile function, the superquantile can be defined as the mean value of the area not covered under the quantile associated with the same probability $\alpha$, that is the right tail of the distribution. Mathematically, the superquantile function is defined as:
\begin{equation}
\label{eq:superquant}
\overline{Q}_{\alpha}(\cdot) = \frac{1}{1-\alpha} \int_{\alpha}^1 Q_{\alpha'}(\cdot) d\alpha' \quad \text{for } \alpha \in (0,1)
\end{equation}

The superquantile formulation can also be adapted to the left tail of the distribution, i.e, $\overline{Q}_{1-\alpha}(\cdot)$. This representation of the CVaR is useful when referring to profits, whereas the right-tail formulation is often used regarding costs. Generally, we will always refer to the superquantile as $\overline{Q}_{\alpha}(\cdot)$, by specifying the side of the tail we are dealing with.

It can be noticed how replacing (\ref{eq:model_intro1_cons1}) with either (\ref{eq:quant}) or (\ref{eq:superq}) is not equivalent, as the employment of the superquantile function gives an extra region of safety, and therefore makes this approach more conservative. However, the superquantile approach, apart from rendering a coherent risk measure, is simpler to implement numerically within optimization problems by using the key result by using the linear formulation proposed by \cite{rockafellar2000optimization}.

As can be seen, this approach has been presented under a stylized context. However, real-world optimization models can present several casuistries that are not directly taken into account in this formulation. For example, from a data-driven perspective, we can benefit from the amount of available information by fitting predictive models to output a direct quantile estimation, getting rid of the conservativeness from the buffered failure probability method. Furthermore, in real complex optimization problems, the decision-maker may not have a direct decision capacity over a certain optimization variable $y$, but there may exist a direct dependency with another decision variable $x$, which can be explicitly modelled. This process is called constraint learning, a topic that is getting attention in the literature lately. For example, in \cite{cremer2018data}, ensembles of decision trees are employed to learn a probabilistic description of a constraint, whereas in \cite{spyros2020decision} stability power system constraints are learned employing neural networks.

In this work, we propose an extension of the constraint learning methodology in order to consider chance constraints, which can also be useful to tackle the inherent uncertainty of point predictive models.

Now suppose we set the optimization model (\ref{eq:model1}). This model aims to minimize an objective cost function $c(x,y|\theta,\xi): \mathbb{R}^{d_x} \rightarrow \mathbb{R}$, which mainly depends on decision variables $x \subset \mathbb{R}^{d_x}$, $y \subset \mathbb{R}^{d_y}$, but at the same time it is conditioned by external contextual information $\theta \subset \mathbb{R}^{d_{\theta}}$ and uncertainty $\xi \subset \mathbb{R}^{d_{\xi}}$.
\begin{subequations}\label{eq:model1}
\begin{align}
\underset{x, y}{\min} \;& \mathbb{E}\left[ c(x,y|\theta;\xi) \right] \\
\text{s.t.}&\notag \\
  &g(x,y|\theta;\xi)\leq 0  \label{eq:model1_cons1}\\
  &y = f^{\mathcal{D}}(x|\theta;\xi) \label{eq:model1_cons2} \\
  &y \leq k \label{eq:model1_cons3}
\end{align}
\end{subequations}

Whereas $x$ represents an accessible, straight-forward decision variable, the nature of the other decision variable $y$ makes it a learned variable. This means that the value set for $y$ is conditioned by the decision we take over $x$ and the contextual information $\theta$. In some contexts, the relationship between $x$, $y$, and $\theta$ can be learned through a selected predictive model $f^{\mathcal{D}}$, which is trained employing a dataset $\mathcal{D} = \{x_i,y_i,\theta_i \}_{i=1}^N$. This process is called ``constraint learning'' \cite{fajemisin2021optimization}, and it is represented by (\ref{eq:model1_cons2}). Furthermore, lets assume that the real-world application requires that an additional restriction has to be imposed over the learned decision variable $y$ in (\ref{eq:model1_cons3}). In this case it is an upper bound constraint, but it can be a lower bound constraint or any other one included within (\ref{eq:model1_cons1})  without loss of generality. These constraints conform the feasible region for $x$ and $y$, i.e., $\mathcal{S}(\theta;\xi)$.

However, the bound constraint (\ref{eq:model1_cons3}) and its fulfillment is directly related to the chosen predictive model $f^{\mathcal{D}}$, how good the training process was, and how accurate the prediction for $y$ is. From an statistical point of view, lets suppose the new prediction for $y$, i.e., $\hat{y}$, follows a estimated Normal distribution $\mathcal{N}(k,\sigma_y)$. In this case, if we employ the point prediction of $\hat{y}$, the optimization problem will consider constraint (\ref{eq:model1_cons3}) fulfilled, as the predictive model will set the value of $y$ equal to $k$. Nevertheless, from a stochastic point of view and taking into account the inherent uncertainty of the predictive model, the bound constraint will not be fulfilled in half of the realizations of $y$ due to its statistical distribution. 

To solve the exposed issues, a Chance Constraint Learning methodology (CCL) for mixed-integer optimization (MIO) problems is developed. The main contributions of this work are four-fold:

\begin{itemize}
    \item[--] to develop a complete data-driven CCL framework, allowing practitioners to train several machine learning models employing their own data, and generate linear constraints to be embedded in their MIO problems.
    \item[--] to extend the standard constraint learning methodology by including chance constraints, and quantile and superquantile estimation methods in order to add robustness to learned decision variables.
    \item[--] to address chance constraints with a direct data-driven functional approach, which avoids using scenario generation and auxiliary binary variables, like SAA-based techniques.
    \item[--] to show the validity of the developed framework in two realistic case studies that verify the good probabilistic performance of the optimal solutions obtained by CCL.
    \item[--] to provide practitioners with an open-access software named ``CCL\_tool'', developed in Python which allows them to implement CCL within their optimization models.
\end{itemize}

The structure of this article is as follows. Section \ref{sec:methodology} describes the proposed methodology in detail, from the conceptual model developed based on linearizable quantile estimation machine learning models to extensions related to superquantiles. Section \ref{sec:software} explains the functioning of the open-access software developed. Sections \ref{sec:case_concrete} and \ref{sec:case_food} test our methodology in two real-world case studies. Finally, Section \ref{sec:conclu} draws the main conclusions of this work.

\section{Chance Constraint Learning Methodology}
\label{sec:methodology}

\subsection{Conceptual model}
\label{sec:conc_model}

We propose an extension of the bound constraint for a learned variable to add probabilistic guarantees. In model (\ref{eq:model2}), the bound constraint (\ref{eq:model1_cons3}) has been replaced by a probabilistic chance constraint (\ref{eq:model2_cons3}), i.e., we seek for the probability of $y$ being lower or equal than $k$ to be greater than a confidence level $\alpha$. 
\begin{subequations}\label{eq:model2}
\begin{align}
\underset{x, y}{\min} \;& \mathbb{E}\left[ c(x|\theta;\xi) \right] \\
\text{s.t.}&\notag \\
  &g(x,y,\theta;\xi)\leq 0  \label{eq:model2_cons1}\\
  &y = f^{\mathcal{D}}(x,\theta;\xi) \label{eq:model2_cons2} \\
  &\mathbb{P}(y \leq k) \geq \alpha \label{eq:model2_cons3}
\end{align}
\end{subequations}

For jointly addressing constraints (\ref{eq:model2_cons2}) and (\ref{eq:model2_cons3}), we propose a simple but effective change over the predictive function $f^{\mathcal{D}}$ in model (\ref{eq:model3}) that will allow to tackle the inherit uncertainty of the point prediction and increase the robustness of the solutions of the optimization problem. The resulting CCL framework is as follows:
\begin{subequations}\label{eq:model3}
\begin{align}
\underset{x, y}{\min} \;& \mathbb{E}\left[ c(x|\theta;\xi) \right] \label{eq:model3_of}\\
\text{s.t.}&\notag \\
  &g(x,y|\theta;\xi)\leq 0  \label{eq:model3_cons1}\\
  &Q_{\alpha}^y = f_{\alpha}^{\mathcal{D}}(x|\theta;\xi) \label{eq:model3_cons2} \\
  &Q_{\alpha}^y \geq k \label{eq:model3_cons3}
\end{align}
\end{subequations}

As can be seen, $f_{\alpha}^{\mathcal{D}}$ can be trained to output an specific conditional quantile of the decision variable $y$: $Q_{\alpha}^y$. Restricting the value of $Q_{\alpha}^y$ to be equal or greater than $k$ makes constraints (\ref{eq:model3_cons2}) and (\ref{eq:model3_cons3}) equivalent to (\ref{eq:model2_cons2}) and (\ref{eq:model2_cons3}). The feasibility of constraint (\ref{eq:model3_cons3}) is then achieved by deciding over decision variable $x$, as it has a direct implication in the value of $y$. Although (\ref{eq:model3_cons2})-(\ref{eq:model3_cons3}) can be viewed as a deterministic counterpart of (\ref{eq:model2_cons2})-(\ref{eq:model2_cons3}), we maintain the expectation in (\ref{eq:model3_of}) to acknowledge that there might be other sources of uncertainty in our model, e.g., uncertain parameters in the objective function and/or in the constraints.

Furthermore, despite conditional quantiles being estimated, this framework can be extended to work with superquantiles, a measure that keeps coherent risk-measures properties. CCL is a flexible framework, that can be employed to tackle uncertainty and give robustness to learned variables by changing point predictive methods by quantile estimation methods, as well as for effective integration in optimization problems where chance constraints are needed due to real-world requirements.

\subsection{Linearizable quantile estimation methods}
\label{sec:mio_models}

Whereas the training of classical point regression methods is based on the minimization of a loss derived from the difference between the point prediction and the actual value of the dependent variable $y$ (mean squared error, the sum of squared errors, etc.), quantile estimation methods are set by minimizing the so-called quantile loss \cite{koenker2017handbook}.

In particular, the $\alpha$-quantile loss $\mathcal{L}_{\alpha}(u)$ is defined in (\ref{eq:quant_loss}), with $u$ representing the residual $y-Q_{\alpha}^y$.
\begin{equation}
\label{eq:quant_loss}
\mathcal{L}_{\alpha}(u) = \begin{cases}
\alpha u, \quad &u \geq 0 \\ (\alpha - 1)u, \quad &u < 0
\end{cases}
\end{equation}

As we can see, this loss is not unique but dependent on the specific estimated $\alpha$-quantile (there will be a loss value for every $\alpha$). When estimating quantiles for a value of $\alpha$ bigger than $0.5$, (\ref{eq:quant_loss}) penalizes the case when the quantile is under $y$ at a greater extent than in the opposite case. On the other side, when $\alpha$ is smaller than $0.5$, being over $y$ brings a bigger penalization.

Employing dataset $\mathcal{D} = \{x_i,y_i,\theta_i \}_{i=1}^N$, the actual value of the $\alpha$-quantile loss $\mathcal{L}_{\alpha}(\mathcal{D})$ is computed by averaging loss values over the complete dataset, as in (\ref{eq:quant_loss_mean}). The minimization of $\mathcal{L}_{\alpha}(\mathcal{D})$ is the main pillar of state-of-art quantile estimation methods. 
\begin{equation}
\label{eq:quant_loss_mean}
\mathcal{L}_{\alpha}(\mathcal{D})  = \frac{1}{N}\sum_{i=1}^{N} \mathcal{L}_\alpha(y_i - Q_{\alpha}^{y_i})
\end{equation}

What follows is an introduction to several quantile estimation methods whose structure is relatively easy to linearize with or without binary variables. This linearization will be critical for embedding these methods as $f_{\alpha}^{\mathcal{D}}(\cdot)$  within mixed-integer linear optimization versions of the problems like (\ref{eq:model3}), and therefore be a participant in a CCL framework.

\subsubsection{Linear Quantile Regression}

One of the simplest approaches is to fit a linear model to estimate the conditional quantile $Q_{\alpha}^{y}$. This model is called Linear Quantile Regression (LQR) and it is based on finding the set of coefficients $\bm{\beta}$ that minimizes the $\alpha$-quantile loss \cite{hao2007quantile}. In this way, LQR outputs $Q_{\alpha}^{y}$ as a function of the decision variables $x$ and contextual information $\theta$.
\begin{equation}\label{eq:lin_quant_reg}
Q_{\alpha}^{y} = \beta_0 + \beta_x x + \beta_{\theta}\theta
\end{equation}

As can be seen, LQR is a method easy to embed within the optimization problem as (\ref{eq:model3_cons2}), allowing also the sensibility and interpretability analysis of the dependent variable with respect to the independent ones.

\subsubsection{Support Vector Quantile Regression}

Support Vector Machine (SVM) is a machine learning method based on hyper-plane cutting in order to generate point predictions or classifications \cite{drucker1996support}. However, the SVM methodology can be adapted by employing the $\alpha$-quantile loss as the minimization target, taking the name of Support Vector Quantile Regression (SVQR). 

We focus on the linear version of SVQR, which makes this model easy to embed in our CCL framework. Linear SVQR is trained in the same way as the exposed LQR, but only penalizes residuals greater than a certain threshold $\epsilon$. After the training is completed, the embedding of the $Q_{\alpha}^{y}$ estimated by the linear SVQR model is done like in LQR (\ref{eq:lin_quant_reg}), i.e., as a linear combination of coefficients $\bm{\beta}$ and the decision variables $x$ and contextual information $\theta$.

\subsubsection{Quantile Regression Tree}
\label{sec:qrt}

A decision tree is a non-linear prediction model based on the creation of sample partitions in different nodes or leaves reached by specific restrictions, and posterior computations to generate predictions. These partitions are made to minimize an error measure. One of the most popular tree-based algorithms is CART (Classification And Regression Tree), firstly introduced in \cite{breiman1984classification}.

Generally, in order to generate a point prediction for a new observation, we follow the tree structure until reaching a certain final leaf. In this leaf, we compute the mean of the dependent variable from the training data samples. This average is assigned to the new observation as a prediction.

However, in the CCL methodology, we focus on predicting conditional quantiles in order to model chance constraints. A simple but effective extension regarding regression trees was proposed in \cite{meinshausen2006quantile}. Basically, instead of computing one average over the data samples in the final leaf, we keep all the observations. With these observations, we can empirically compute every estimated conditional quantile for the dependent variable, i.e., $Q_{\alpha}^{y}$.

\begin{figure}%[ht]
    \centering
    \includegraphics[width=\textwidth]{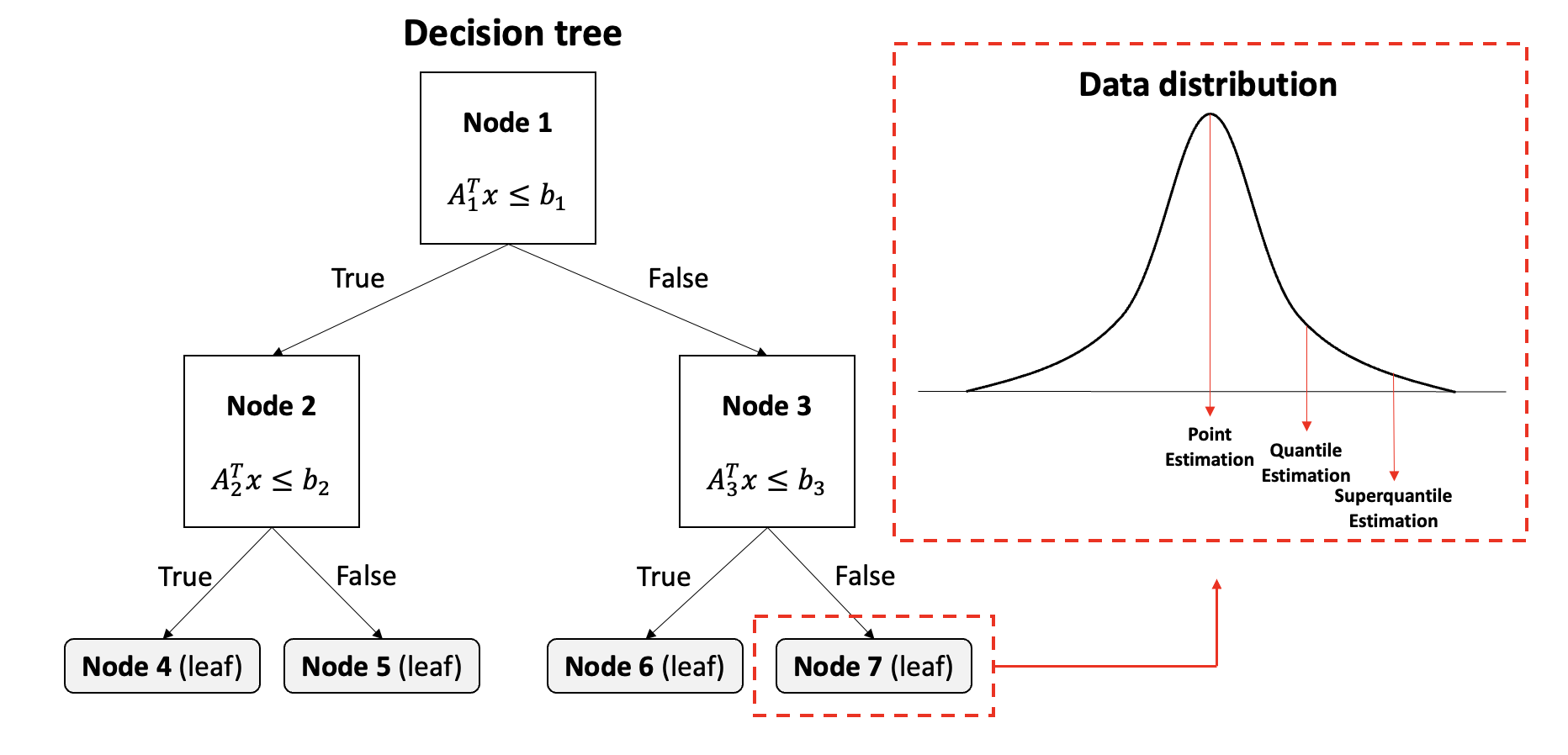}
    \caption{Decision tree regression structure.}
    \label{fig:tree_str}
\end{figure}

Figure \ref{fig:tree_str} shows the typical structure of a decision tree (in this case a tree of depth 2). As can be seen, to make a prediction, we move down the tree following the path created by some restrictions. When we arrive at a final node (a leaf), we can find the training data distribution. For example, for an observation $x$ to be assigned to the fourth node, the set of inequalities $\mathcal{N}_4: \{A_1^T x \leq b_1,  A_2^T x \leq b_2\}$ must be true. In point prediction, we would compute an average to get the expected value of the prediction when a leaf node is reached. Nevertheless, under the CCL framework, we can actually compute the empirical quantile, or even the superquantile as the average of the samples greater than the quantile (or below it without loss of generality).

In order to embed a tree-based prediction model in a mixed-integer optimization problem, we have to write these restriction linearly to obtain an equivalent set of constraints. For that purpose, we follow the formulation proposed in \cite{maragno2021mixed} with a slight modification to include quantile estimation. Thus, considering the tree in Figure \ref{fig:tree_str}, the leaf assignment encoding for an observation $x$ is described as: 
\begin{subequations}\label{eq:tree_enco}
\begin{align}
&A_1^T x - M(1-u_4) \leq b_1  \label{eq:tree_cons1} \\
&A_2^T x - M(1-u_4) \leq b_2  \label{eq:tree_cons2} \\
&A_1^T x - M(1-u_5) \leq b_1  \label{eq:tree_cons3} \\
-&A_2^T x - M(1-u_5) \leq -b_2 - \epsilon  \label{eq:tree_cons4} \\
-&A_1^T x - M(1-u_6) \leq -b_1 - \epsilon \label{eq:tree_cons5} \\
&A_3^T x - M(1-u_6) \leq b_3  \label{eq:tree_cons6} \\
-&A_1^T x - M(1-u_7) \leq -b_1 - \epsilon \label{eq:tree_cons7} \\
-&A_3^T x - M(1-u_7) \leq -b_3 - \epsilon  \label{eq:tree_cons8} \\
&u_4 + u_5 + u_6 + u_7 = 1  \label{eq:tree_cons9} \\
&y - (q_4^{\alpha}u_4 + q_5^{\alpha}u_5 + q_6^{\alpha}u_6 + q_7^{\alpha}u_7) = 0  \label{eq:tree_cons10}
\end{align}
\end{subequations}

\noindent where $u_i$ represent a binary variable for the corresponding leaf $i$, $\epsilon$ a small parameter to remove strict inequalities, and $M$ a sufficiently large constant (big-M approach).

For $x$, if $A_1^T x \leq b_1$, constraints (\ref{eq:tree_cons5}) and (\ref{eq:tree_cons7}) will force $u_6$ and $u_7$ to be zero. Furthermore, if $A_2^T x \leq b_2$, constraint (\ref{eq:tree_cons4}) will assign $u_5$ a value of zero. Finally, constraint (\ref{eq:tree_cons9}) will force $u_4$ to be one, assigning to $y$ (\ref{eq:tree_cons10}) the value of the predicted quantile $q_4^{\alpha}$ (the quantile in leaf 4).

\subsubsection{Tree-based quantile ensemble methods}

An ensemble is a machine learning method based on generating predictions from multiple independent base simple models. In this section, we will focus on two ensemble methods based on quantile regression trees: Quantile Regression Forest (QRF) and Gradient Boosting Quantile Regression (GBQR).

QRF \cite{meinshausen2006quantile} is a method closely related to randomization, as data and features are bootstrapped to grow the quantile regression tree that can be chosen also randomly. Typically, QRF generates a new prediction by averaging over the prediction of each base model, i.e., $Q_{\alpha}^{y_i} = \frac{1}{T}\sum_{t\in T} Q_{\alpha, t}^{y_i} $, where $T$ represents the number of quantile regression trees conforming the ensemble. 

The functioning of each tree is the same as the one exposed in Section \ref{sec:qrt}. Therefore, in order to embed QRF in a MIO problem, we have to embed each tree one by one, while adding a final constraint to generate a final quantile estimation as the average of all the estimations from the base models.

GBQR is an extension from the gradient boosting (GB) methodology proposed in \cite{friedman2002stochastic}. GB aims to improve the performance of a weak base model (``to boost'') by training a new ensemble member with a modified dataset, i.e., with the residuals from the previous base model instead of the original dependent variable. GB can be implemented with any regression method as a base learner, although it is mostly used within regression trees in the literature. To achieve a conditional quantile estimation, tree-based models of the ensemble are grown making use of the quantile loss (\ref{eq:quant_loss}), taking the name of GBQR.

\subsubsection{Quantile Regression Neural Networks}

Neural Networks (NNs) are powerful and flexible techniques employed to a great extent in the current machine learning and artificial intelligent state-of-art \cite{gallant1993neural}. NNs are composed of one input layer, $L-2$ hidden layers, and one final output layer. Each layer can be composed of several nodes or neurons that perform computations regarding weights and biases.

In a given hidden layer $l$ of the NN, with nodes $N^l$, the value of a node $i \in N^l$, denoted as $v_i^l$ , is calculated using the weighted sum of the previous layer’s node values, followed by a non-linear activation function $g(\cdot)$. This value is given as:
\begin{equation*}
\label{eq:output_1}  
 v_i^{l}=g\left(b_i^l + \sum_{j\in N^{l-1}} w_{ij}^{l} v_j^{l-1} \right)
\end{equation*}

\noindent where $b_i^l$ represents the bias term, and $w_{ij}^{l}$ the weights or coefficients for node $i$ in layer $l$. This process will continue in the following layers, taking as inputs the outputs of the previous layers, until reaching the final one.

Recent works \cite{grimstad2019relu, anderson2020strong} have studied  rectified linear unit (ReLU)-based NNs as a MIO-representable class of neural networks, and therefore, adequate for the CCL framework. Taking the previous example, if $g(\cdot)$ represents the ReLU activation function, the output from node $i$ in layer $l$ is computed as:
\begin{equation*}
\label{eq:output_2}  
 v_i^{l}=\max\left(0, b_i^l + \sum_{j\in N^{l-1}} w_{ij}^{l} v_j^{l-1} \right)
\end{equation*}

The training process of this NN can be focused on minimizing the $\alpha$-quantile loss in order to output the corresponding conditional quantile $Q_{\alpha}^{y_i}$. This approach is denoted as Quantile Regression Neural Networks (QRNN). Thus, constraints for the CCL framework can be generated in a recursive procedure from the input layer to the output layer. This final layer $L$ will consist in an unique neuron which, by means of a linear combination of values obtained in layer $L-1$, will output $Q_{\alpha}^{y}$ as follows:
\begin{equation*}
\label{eq:output_3}  
Q_{\alpha}^{y} =  b^L + \sum_{j\in N^{L-1}} w_{j}^{L} v_j^{L-1} 
\end{equation*}

We can anticipate that the complexity and the number of constraints will grow exponentially using QRNNs instead of linear models such as LQR or linear SVQR. Furthermore, as the number of layers increases, the explainability of the model will decrease, transforming QRNN into a ``black-box'' model. However, if we focus on prediction accuracy, the use of NNs gives a better performance than the use of linear models in many relevant applications.

\subsection{Empirical superquantile estimation}
\label{sec:superq_estim}

As it has been mentioned in Section \ref{sec:intro}, the main advantage of using estimated empirical quantiles (VaR) to deal with probabilistic constraints is getting rid of the extra buffered probability associated with the superquantile (CVaR). However, CVaR possesses some modelling advantages as it is a coherent risk measure.

Thus, we propose an extension of the CCL methodology in order to employ superquantiles, so that the decision-maker has the option to choose between VaR and CVaR regarding its preferences for risk modelling. 

This derivation is based on Mixed Quantile Regression (\ref{eq:mixed_qr}), a methodology proposed in \cite{rockafellar2008risk} and \cite{chun2012conditional}. In particular, for linear problems:
\begin{equation}
\label{eq:mixed_qr}
\overline{Q}_{\alpha}^{y_i} = \text{CVaR}_{\alpha}[y|x,\theta] = \bm{\overline{\beta}}(x,\theta)
\end{equation}

This implies that for an empirical approximation of the superquantile we can use $\sum_{j=1}^M w_j Q_{\alpha_j}^{y_i}$ which is based on the discretization of the probability distribution of $y$ employing $M$ different quantiles (or VaRs). Results from \cite{harsha2021prescriptive} show that, for heteroscedastic datasets, $\bm{\overline{\beta}} = \sum_{j=1}^M w_j \bm{\beta}_j$. This means that coefficients for a CVaR-based model can be obtained by means of a linear combination of coefficients from $M$ different linear quantile models. Extensions for homocedastic datasets can be found in the same work. Moreover, to determine the values of weights $w_j$ as well as the quantiles employed for modelling the CVaR, a midpoint quadrature rule can be used. In this way, with $\Delta = J^{-1}(1-\alpha)$, we have $w_j = (1-\alpha)^{-1}\Delta$ and $\alpha_j = \alpha + (j-0.5)\Delta, \; j=1,\dots,J$ as the weighs and VaRs to employ, when dealing with the right tail of the distribution. In the case of a CVaR referring to the left side of the distribution, we can use $\Delta = J^{-1}\alpha$, and $w_j = \alpha^{-1}\Delta$ and $\alpha_j = \alpha - (j-0.5)\Delta, \; j=1,\dots,J$ as the weighs and quantiles, respectively.

However, the above procedure can not be extended to the case of tree-based prediction models. In this case, there are no coefficients $\bm{\beta}$ available, but data partitions in leaf nodes. For that reason, we propose to empirically compute the superquantile from the data distribution presented in leaf nodes. Taking as an example Figure \ref{fig:tree_str}, we can see how, from the data distribution of leaf node 7, the empirical quantile can be computed. From that, the superquantile is calculated as the mean value from the observations greater than the quantile (in the case of a right tail CVaR). To give robustness to the estimation, there must be enough observations in the leaf. Therefore, we should control the minimum number of samples within leaves as a hyper-parameter of our tree-based models. This approximation can be employed for decision trees and ensembles such as QRF.

Finally, for QRNNs, although there are coefficients $\bm{\beta}$ associated with each neuron, the methodology proposed in \cite{harsha2021prescriptive} cannot be applied when we include several layers and neurons, as we turn the model into a non-linear one. Thus, we propose to directly compute the conditional superquantile as the weighted sum of estimated quantiles, i.e., $\overline{Q}_{\alpha}^{y} = \sum_{j=1}^M w_j Q_{\alpha_j}^{y}$. For instance, we can make use of the same weights values and quantiles of the midpoint quadrature rule. Furthermore, we can accelerate the training process and fit one unique QRNN that outputs the $J$ seek quantiles in the final layer, instead of one NN for each quantile. In this way, all the quantiles from the NN will share all of the coefficients except those from the final layer. However, obtaining multiple quantile outputs in a single NN can bring up the well-known issue of quantile crossing, that is, the resulting distribution function may not be monotonically increasing. To ensure that the resulting QRNN solves this issue, we add a penalty term in case of quantile crossing to the quantile loss during the training phase \cite{moon2021learning}.

As can be seen, there is not a well-defined homogeneous methodology for computing the superquantile regardless of the prediction model. However, we believe that the proposed methodologies that numerically estimate the CVaR can be of great interest to practitioners working on real-world problems. These proposals are validated empirically in Sections \ref{sec:case_concrete} and \ref{sec:case_food}.

\section{CCL software development}
\label{sec:software}

In this section, we introduce the software developed in order to implement the chance constraint learning methodology.

The ``CCL\_tool'' has been fully developed in Python, and it is completely open access for practitioners in \cite{Alcantara_ccl_tool}. This tool is Scikit-learn \cite{scikit-learn} and Pytorch \cite{NEURIPS2019_9015} based for fitting different machine learning models, while also extending some code from \cite{OptiCL} for embedding predictive models within MIO problems modeled in Pyomo \cite{hart2017pyomo}.

\begin{figure}[!ht]
    \centering
    \includegraphics[width=\textwidth]{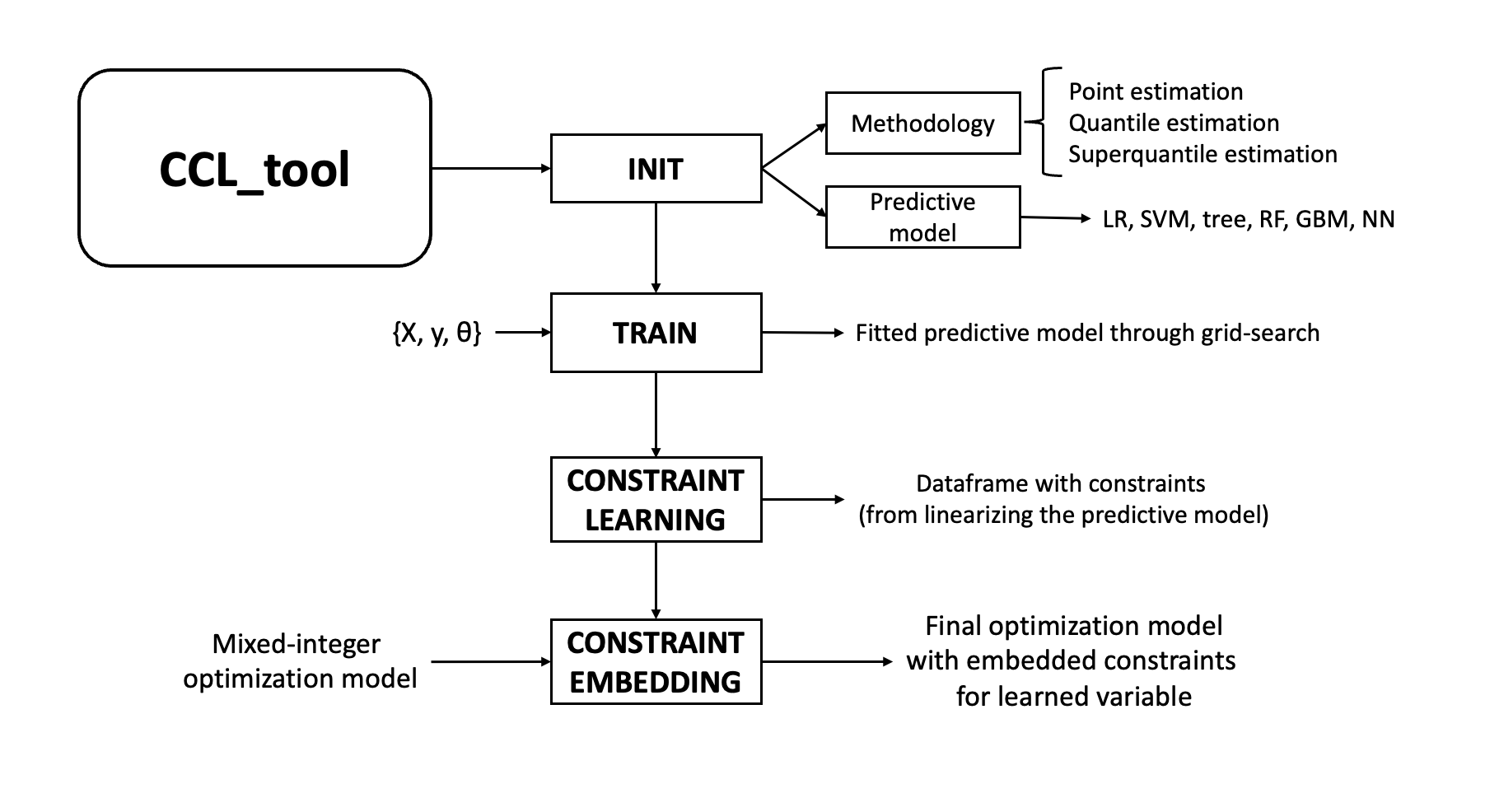}
    \caption{Chance constraint learning tool structure}
    \label{fig:ccl_tool}
\end{figure}

Its functioning is defined by four steps. First, an initialization where the learning methodology and prediction model are selected by the user. Then, a training process where, from a dataset, we fit the predictive model. Later, once the model is trained, we generate the learned constraints from the predictive model and finally, these constraints are embedded in a MIO problem. A summary of this structure is sketched in Figure \ref{fig:ccl_tool}.

In the initialization phase, we choose a methodology between point, quantile, or superquantile estimation. That is, with point estimation we want a model to estimate the expected value of the learned variable (with this variable having a possible bound). However, if quantile or superquantile estimation is chosen, we define a CCL problem, i.e., we want our learned variable to be greater or lower than a value with a certain probability. In this case, we will have to select the confidence level $\alpha$, and the number of quantiles to approximate the superquantile depending on the method (see Section \ref{sec:superq_estim}). Regarding the predictive models, linear regression, tree models, random forest, and (deep) neural networks are available for both point, quantile, and superquantile estimation. Only GB is not available for superquantile estimation as it is being developed from scratch, and SVM for (super)quantile estimation, as there is no Python-based package for SVM quantile estimation with linear kernel. We suggest practitioners interested in SVM train models in R and embed the model as linear regression. 

When we enter the training step, a dataset $\mathcal{D} = \{x_i,y_i,\theta_i\}_{i=1}^N$ is required, with $x$ representing past decisions, $y$ being the learned variable, and $\theta$ the contextual information. In this way, the tool will split the dataset into a training and validation subset to select the best hyper-parameter setting for the predictive model. For example, the number of trees, the maximum depth, or the minimum number of samples in leaf for Random Forest. This selection is the one that minimizes the validation error, for example, the MAE for point estimation and the $\alpha$-quantile loss (\ref{eq:quant_loss_mean}) for quantile and superquantile estimation.

Once the model is fitted, constraints are generated as exposed in Section \ref{sec:mio_models} depending on the selection of the predictive model. These constraints are finally added to our MIO problem.

Therefore, a fully data-driven CCL methodology is developed, where practitioners can find optimal solutions to their decision-making problems under uncertainty by using high-level machine learning models and information while solving chance constraints and hedging the risk associated with the value of a non-accessible decision variable: the learned constraint.

\section{Case Studies}
\label{sec:cases}

In this section, we empirically validate the proposed CCL methodology in two realistic frameworks. For clarity, in these case studies, no contextual information $\theta$ is employed to build the learned constraints. However, we can easily deal with this extension in the ``CCL\_tool'', in particular, by setting $\theta$ as decision variables within vector $x$, and then by fixing their values in the MIO problem to the observed ones.

\subsection{Case study I: concrete compressive strength}
\label{sec:case_concrete}

In this case study, we develop an optimization problem based on the very well know ``concrete'' dataset, first time published in \cite{yeh1998modeling}. This dataset contains information about concrete compressive strength (dependent variable) measured in megapascals (MPa), along with values from the concrete components and age (independent and decision variables).

Establishing the relationship between the compressive strength and its components through a linearizable prediction model will allow us to embed this learned variable into mixed-integer optimization problems like the one we propose. During the development of this case study, we will focus on the different realizations the decision variables take depending on the constraint learning methodology. That is, we study the decision-maker behaviour from the (point) constraint learning and the chance constraint learning (using both quantiles and superquantiles) perspectives. 

\subsubsection{Optimization problem}

The optimization problem will be based on minimizing concrete production costs while ensuring a predefined compressive strength level. We assume the producer has complete flexibility to decide when to produce, and the production cost to be measured in $kg/m^3$, to preserve the nature of the dataset, as it will be detailed later.

First, we define the following notation for the problem:

\medskip
\noindent Indices and sets:

\begin{itemize}
    \item[--] $I$: Set of concrete components, indexed by $i$.
\end{itemize}

\medskip
\noindent Variables:
\begin{itemize}
	\item[--] $x_{i}$: Quantity of concrete component $i$, measured in kilograms, employed for producing one $m^3$.
	\item[--] $d$: Age of strength testing (in days).
	\item[--] $y_{\text{strg}}$: Concrete compressive strength, measured in MPa.
\end{itemize}

\medskip
\noindent Parameters:
\begin{itemize}
	\item[--] $c_i$: Cost in euros for concrete component $i$, per kilogram.
	\item[--] $\underbar{\text{weight}}$ and $\overline{\text{weight}}$: Lower and upper bound for the total component weight used in the concrete production, respectively.
	\item[--] $k$: Concrete compressive strength lower bound.
	\item[--] $\alpha$: Probability confidence level.
\end{itemize}

The formulation of the chance constrained optimization problem is presented below:
\begin{subequations}\label{eq:model_concrete}
\begin{align}
\underset{\bm{x}, y_{\text{strg}}, d}{\min} \;& \sum_{i\in I} c_i x_{i}   \label{eq:model_concrete_ob}\\
\text{s.t.}&\notag \\
  &y_{\text{strg}} = f^{\mathcal{D}}(\bm{x},d) \label{eq:model_concrete_cons1} \\
  &\mathbb{P}(y_{\text{strg}} \geq k) \geq \alpha  \label{eq:model_concrete_cons2} \\
  &\underbar{\text{weight}} \leq \sum_{i\in I} x_{i} \leq \overline{\text{weight}}  \label{eq:model_concrete_cons3}
\end{align}
\end{subequations}

The objective function of this problem (\ref{eq:model_concrete_ob}) aims to minimize the production costs. The producer will have to decide how much quantity of each component $x_{i}$ add to the mixture in order to produce one $m^3$ of concrete. Furthermore, he must decide when to produce, as the compressive strength is dependent on the age of the concrete.

Constraint (\ref{eq:model_concrete_cons1}) is the so-called ``learned constraint'', which defines the ``learned decision variable'': the compressive strength for a specific concrete mixture, $y_{\text{strg}}$. We can see how this strength is dependent of the components $\bm{x}$ and the age of the concrete $d$. This relationship will be established through a predictive function $f^{\mathcal{D}}$, trained using a dataset $\mathcal{D}$, and embedded within the optimization problem. Next, constraint (\ref{eq:model_concrete_cons2}) defines a chance constraint. The producer wants its concrete to have a compressive strength $y_{\text{strg}}$ at least equal to a defined value $k$, with a certain probability $\alpha$. Notice how both constraints (\ref{eq:model_concrete_cons1}) and (\ref{eq:model_concrete_cons2}) create a CCL framework together.

Finally, constraint (\ref{eq:model_concrete_cons3}) ensures that the weight of the concrete components is between a lower and upper bound.

\subsubsection{Methodology}

As it has been mentioned before, the well-known ``concrete'' dataset will be employed in this problem. Table \ref{tab:vars_concrete} indicates the variables presented in this dataset composed by $1030$ samples. We can see how all the components are measured in $kg/m^3$, so that we assume orders are made in $m^3$.

\begin{table}[!ht]
    \centering
    \caption{Set of variables employed for the concrete compressive strength problem.}
    \label{tab:vars_concrete}
    \begin{tabular}{lll}
        \multicolumn{1}{c}{\textbf{Type of variable}} & \multicolumn{1}{c}{\textbf{Name}} & \multicolumn{1}{c}{\textbf{Units}}  \\ \hline
        Dependent & Compressive strength & Mpa  \\ \hline
        \multirow{8}{*}{Independent}  & Cement & $kg/m^3$  \\ 
         & Blast Furnace Slag & $kg/m^3$  \\
         & Fly Ash & $kg/m^3$  \\
         & Water & $kg/m^3$  \\ 
         & Superplasticizer & $kg/m^3$  \\ 
         & Coarse Aggregate & $kg/m^3$  \\ 
         & Fine Aggregate & $kg/m^3$  \\ 
         & Age & Days \\ \hline
    \end{tabular}
\end{table}

The next information needed for the optimization problem are the costs associated with each concrete component. This information is not provided in the original dataset, so we assume the producer to have the costs presented in Table \ref{tab:cost_concrete}.

\begin{table}[!ht]
    \centering
    \caption{Concrete components costs}
    \label{tab:cost_concrete}
    \begin{tabular}{lc}
        \textbf{Component} & \textbf{Cost} [$\text{\euro}$\textbf{/kg}] \\ \hline
        Cement & $0.050$  \\ 
        Blast Furnace Slag & $0.040$  \\
        Fly Ash & $0.045$  \\
        Water & $0.002$  \\ 
        Superplasticizer & $1.800$  \\ 
        Coarse Aggregate & $0.020$  \\ 
        Fine Aggregate & $0.020$  \\ \hline
    \end{tabular}
\end{table}

Regarding the predictive task to illustrate the CCL methodology, Random Forest (RF) will be employed as the predictive model to relate compressive strength and concrete components. Concerning the different methodologies, point estimation will be employed to ensure a compressive strength bigger than $45$ Mpa without taking into account the probability of chance constraints. On the other side, quantile and superquantile estimation methodology will be employed under the CCL framework to ensure compressive strength to be bigger than $45$ Mpa with a probability bigger than $95\%$. This means that, as a lower bound restriction is presented in the problem, we need to model a quantile and superquantile at 5\% $(1-\alpha)$ level and force it to be greater than $45$ Mpa. Then, the RF will estimate the conditional quantile (or superquantile) at level 5\%.

Finally, the lower and upper bounds, $\underbar{\text{weight}}$ and $\overline{\text{weight}}$ for the component weights in the mixture will be $2230$ kg and $2450$ kg, respectively (chosen as the fifth and ninety-fifth quantile of the component weight in the dataset).

\subsubsection{Results}
\label{sec:concrete_results}

The optimal concrete components mixture problem has been solved through a Python 3.9.12 implementation, using Pyomo 6.3 \cite{hart2017pyomo}. The selected mathematical solver for all the computations was Gurobi \cite{gurobi} in its version 9.5. Besides, the computer employed included a CPU Intel Core i7 10700, RAM of 64 GB, and NVIDIA GeForce GTX 2060 graphic card. As we mentioned before, the optimal solution will be obtained for the three different exposed methodologies: point estimation, quantile estimation, and superquantile estimation methodology.

Table \ref{tab:concrete_sol} shows the optimal concrete mixture for ensuring a compressive strength depending on the selected methodology when choosing RF as the predictive model. As can be seen, the cement quantity increases as we move to more restrictive methodologies to ensure compressive strength. The opposite occurs with the quantity of water and coarse aggregate. It is interesting to notice the null use of superplasticizer (a really expensive component) and the differences in the age of the concrete for the different methodologies.

\begin{table}[!ht]
    \centering
    \caption{Optimal concrete components mixture}
    \label{tab:concrete_sol}
    \begin{tabular}{c|ccc}
        \multicolumn{1}{p{3cm}|}{\centering \textbf{Decision Variable} } & \multicolumn{1}{p{3cm}}{\centering \textbf{Point} \\ \textbf{Estimation} } & \multicolumn{1}{p{3cm}}{\centering \textbf{Quantile} \\ \textbf{Estimation} } & \multicolumn{1}{p{3cm}}{\centering \textbf{Superquantile} \\ \textbf{Estimation} } \\ \hline
        Cement & 166.90 kg & 357.50 kg & 357.50 kg  \\ 
        Blast Furnace Slag & 43.72 kg & - & 59.35 kg  \\
        Fly Ash & - & - & -  \\
        Water & 175.08 kg & 176.45 kg & 153.7 kg \\ 
        Superplasticizer & - & - & - \\ 
        Coarse Aggregate & 1149.93 kg & 1016.7 kg & 945.75 kg \\ 
        Fine Aggregate & 694.37 kg & 679.35 kg & 713.7 kg  \\ 
        Age & 78 days & 95 days & 73 days \\ \hline
        \textbf{Cost} & 47.33 \text{\euro} & 52.15 \text{\euro} & 53.75 \text{\euro} \\ \hline
    \end{tabular}
\end{table}

The biggest employment of cement in the quantile and superquantile methodologies suggests that the optimal solution would be more expensive, i.e., we need to ensure the compressive strength to a greater extent. For that reason, we analyze the objective value for the three different cases. For the point estimation methodology, a cost of 47.33 \text{\euro} for the concrete mixture is obtained. In the case of quantile estimation, the cost rises to 52.15 \text{\euro}. Finally, the biggest cost is obtained by making use of the superquantile methodology, with a value of 53.75 \text{\euro}. We can see an increase in the cost of more than 6 \text{\euro} (almost 15\%) between the point and the superquantile estimation. This difference in the cost should imply an increase in the safety regarding the compressive strength to support the (super)quantile estimation approach. To ensure that, we compute the distribution function of the three different optimal solutions. In that sense, we fit a Quantile Regression Forest with our training dataset. With this fitted model, we will predict the conditional quantile of the optimal solution at several values of $\alpha$ to recreate the distribution function.

\begin{figure}[!ht]
    \centering
    \includegraphics[width=\textwidth]{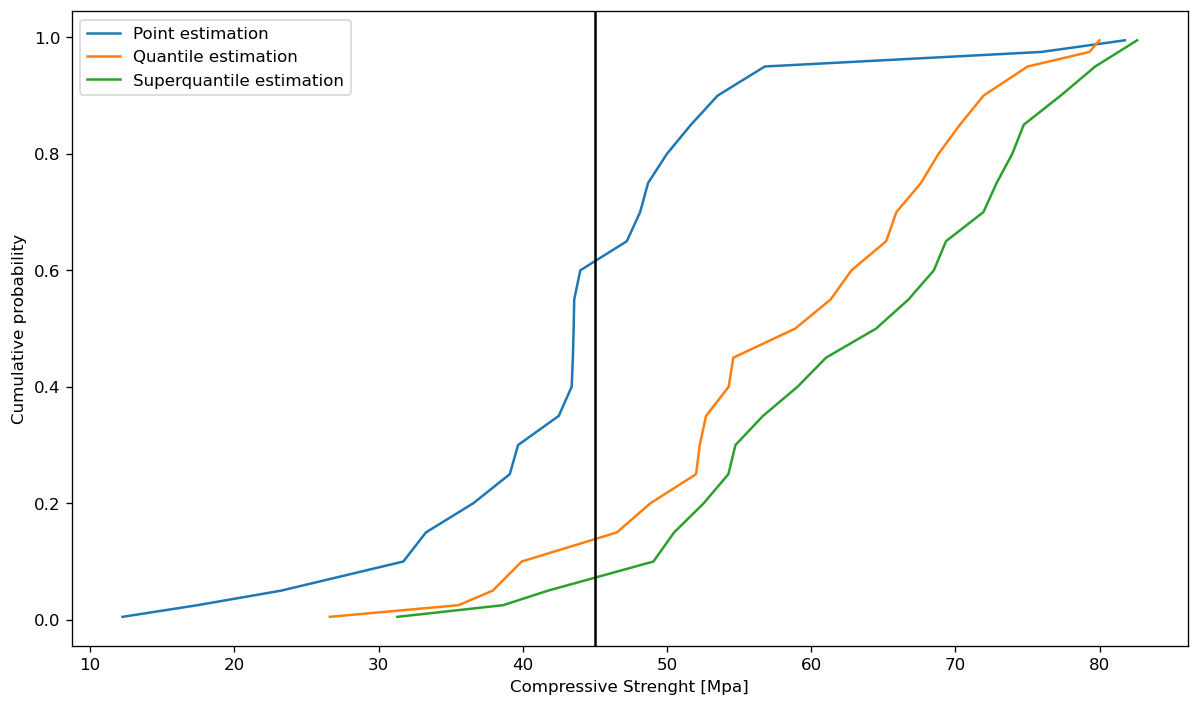}
    \caption{Optimal solution's distribution function for the compressive strength problem regarding the employed methodology.}
    \label{fig:concrete_sol_dist}
\end{figure}

In Figure \ref{fig:concrete_sol_dist}, blue, orange, and green curves represent the distribution function for the optimal solution obtained from the point estimation, quantile estimation, and superquantile estimation methodology, respectively. The vertical black line indicates the value of the required compressive strength in the optimization problem (45 Mpa).

We can see how, for point estimation methodology, most of the probability distribution is under the required compressive strength. This shows how, even fitting a powerful machine learning point prediction model, the distribution of the optimal solution can be biased. Note that, in most of the real-life scenarios, the actual compressive strength would be below 45 Mpa.

On the other side, the distribution of the optimal solution for the quantile and superquantile methodology is quite shifted to the right with respect to the one from the point methodology. In these cases, most of the distribution is over the required compressive strength.

Before concluding this case study, if we go beyond RF, Table \ref{tab:concrete_perf} shows a brief performance comparison between different predictive models (Section \ref{sec:mio_models}) regarding their testing error, fitting, and constraint learning time, optimization problem resolution time, number of constraints, binary variables, and continuous variables, within the final optimization problem. For RF models, a grid search methodology is applied to select the best set of hyper-parameters, whereas, for NN predictive models, a 100 neurons single-layer structure is selected.

\begin{table}[!ht]
    \centering
    \caption{Chance constraint learning performance regarding methodology and predictive model for the compressive strength problem}
    \label{tab:concrete_perf}
    \resizebox{\textwidth}{!}{%
    \begin{tabular}{l|c|cccccc}
        \multicolumn{1}{p{3cm}|}{\centering \textbf{Methodology} } & \multicolumn{1}{p{2cm}|}{\centering \textbf{Predictive model} } & \multicolumn{1}{p{4cm}}{\centering \textbf{Testing error (MAE/Quant loss)} } & \multicolumn{1}{p{2.5cm}}{\centering \textbf{Fitting \& CL time (s)} } & \multicolumn{1}{p{3.5cm}}{\centering \textbf{Optimization solving time (s)} } & \multicolumn{1}{p{3cm}}{\centering \textbf{Problem size (constraints)} } & \multicolumn{1}{p{2cm}}{\centering \textbf{Binary vars} } & \multicolumn{1}{p{2.25cm}}{\centering \textbf{Non-binary vars}} \\ \hline
        \multirowcell{3}{Point\\Estimation}  & LR & 7.87 & 0.5 & 0.3 & 4 & 0 & 9 \\
        ~ & RF & 4.86 & 61.5 & 274.2 & 454 & 3714 & 159 \\ 
        ~ & NN & 3.41 & 58.9 & 0.3 & 304 & 100 & 109 \\ \hline
        \multirowcell{3}{Quantile\\Estimation} & LR & 0.94 & 0.1 & 0.2 & 4 & 0 & 9\\ 
        ~ & RF & 0.70 & 105.8 & 84.4 & 604 & 3711 & 159\\
        ~ & NN & 0.42 & 63.1 & 0.3 & 304 & 100 & 109\\ \hline
        \multirowcell{3}{Superquantile\\Estimation} & LR & 0.96 & 0.7 & 0.2 & 4 & 0 & 9 \\
        ~ & RF & 0.70 & 109.2 & 50.4 & 604 & 3711 & 159\\ 
        ~ & NN & 0.32 & 76.7 & 0.3 & 309 & 100 & 114\\ \hline
    \end{tabular}%
    }
\end{table}

First of all, we can see how linear and NN-based methods are the fastest in obtaining the optimal solution. Besides, we can notice the higher number of constraints and variables we have to deal with when complex models like RF and NN are employed as our predictive models. Furthermore, due to the high number of binary variables that RF models include within the optimization problem, the optimization solving time when using these tree-based methods is bigger to a higher extent. As an additional note, we cannot compare the fitting process time between RF and NN, as hundreds of RF are trained during the grid search phase. Generally, the NN-based model obtains the lowest error (MAE or quantile loss) regarding the test partition of our dataset. However, the small number of samples in the dataset makes NNs not generalize appropriately (overfitting).

\subsection{Case study II: palatable food basket}
\label{sec:case_food}

In this section, we address the case study presented in \cite{maragno2021mixed}; being a simplification of the one presented in \cite{peters2021nutritious}. In particular, a decision-making humanitarian operation tool is provided to the World Food Program (WFP) to optimize, among others, the food basket to be delivered. A food basket has to ensure certain nutrition requirements while also accomplishing some palatability rules. In \cite{maragno2021mixed}, palatability scores are set to be learned by a constraint learning process, employing machine learning models for point estimation. In this work, we propose to include chance constraints for the model to ensure a minimum level of palatability with a certain probability.

Therefore, we will obtain and compare three different optimal food baskets regarding the three addressed methodologies of this work: point, quantile, and superquantile estimation.

\subsubsection{Optimization problem}

The optimization problem is based on a network flow model, restricted by nutrition requirements and food basket palatability. In this problem, we aim to minimize the procurement and transportation costs of the different commodities.

In the following, the notation for the problem is introduced.

\medskip
\noindent Indices and sets:

\begin{itemize}
    \item[--] $\mathcal{N_S}$: Set of source nodes.
    \item[--] $\mathcal{N_T}$: Set of transshipment nodes.
    \item[--] $\mathcal{N_D}$: Set of delivery nodes.
    \item[--] $\mathcal{K}$: Set of commodities, indexed by $k$.
    \item[--] $\mathcal{L}$: Set of nutrients, indexed by $l$.
\end{itemize}

\medskip
\noindent Variables:
\begin{itemize}
	\item[--] $F_{ijk}$: Metric tons of commodity $k$ transported between node $i$ and node $j$.
	\item[--] $x_k$: Grams of commodity $k$ in the food basket.
	\item[--] $y_{\text{pltb}}$: Food basket palatability
\end{itemize}

\medskip
\noindent Parameters:
\begin{itemize}
	\item[--] $\gamma$: Conversion rate from metric tons (mt) to grams (g).
	\item[--] $d_i$: Number of beneficiaries at delivery point $i \in \mathcal{N_D}$
	\item[--] $days$: Number of feeding days.
	\item[--] $Nutreq_l$: Nutritional requirement for nutrient $l \in \mathcal{L}$ (grams/person/day).
	\item[--] $Nutval_{kl}$: Nutritional value for nutrient $l \in \mathcal{L}$ per gram of commodity $k \in \mathcal{K}$.
	\item[--] $p^P_{ik}$: Procurement cost (in \$/mt) of commodity $k$ from source $i \in \mathcal{N_S}$.
	\item[--] $p^T_{ijk}$: Transportation cost (in \$/mt) of commodity $k$ from node $i \in \mathcal{N_S} \cup \mathcal{N_T}$ to node $j \in \mathcal{N_T} \cup \mathcal{N_D}$.
	\item[--] $t$: Palatability score lower bound.
	\item[--] $\alpha$: Probability confidence level.
\end{itemize}

The extension of the palatable food basket problem to include chance constraints is presented as follows:
\begin{subequations}\label{eq:model_food}
\begin{align}
\underset{\bm{x}, y_{\text{pltb}}, \bm{F}}{\min} \;& \sum_{i\in \mathcal{N_S}} \sum_{j \in \mathcal{N_T} \cup \mathcal{N_D}} \sum_{k \in \mathcal{K}} p^P_{ik}F_{ijk} + \sum_{i\in \mathcal{N_S} \cup \mathcal{N_T}} \sum_{j \in \mathcal{N_T} \cup \mathcal{N_D}} \sum_{k \in \mathcal{K}} p^T_{ijk}F_{ijk}  \label{eq:model_food_ob}\\
\text{s.t.}&\notag \\
  &\sum_{j \in \mathcal{N_T}} F_{ijk} = \sum_{j \in \mathcal{N_T}} F_{jik}, \quad i \in \mathcal{N_T}, k \in \mathcal{K} \label{eq:model_food_cons1}\\
  &\sum_{i\in \mathcal{N_S} \cup \mathcal{N_T}} \gamma F{jik} = d_ix_k days, \quad i \in \mathcal{N_D}, k \in \mathcal{K} \label{eq:model_food_cons2}\\
  &\sum_{k \in \mathcal{K}} Nutval_{kl} x_k \geq Nutreq_l, \quad l \in \mathcal{L} \label{eq:model_food_cons3}\\
  &x_{salt} = 5 \label{eq:model_food_cons4}\\
  &x_{sugar} = 20 \label{eq:model_food_cons5}\\
  &y_{\text{pltb}} = f^{\mathcal{D}}(\bm{x}) \label{eq:model_food_cons6} \\
  &\mathbb{P}(y_{\text{pltb}} \geq t) \geq \alpha  \label{eq:model_food_cons7}
\end{align}
\end{subequations}

The objective function (\ref{eq:model_food_ob}) aims to minimize both procurement and transportation costs. Constraint (\ref{eq:model_food_cons1}) is employed to balance the network. Then, constraint (\ref{eq:model_food_cons2}) ensures that the flow into the delivery node is equal to the demand, computed as the number of beneficiaries times the daily quantity of commodity $k$ times the number of days. Constraint (\ref{eq:model_food_cons3}) forces the optimal food basket to meet the nutritional requirements. Constraints (\ref{eq:model_food_cons4}) and (\ref{eq:model_food_cons5}) maintain the amount of salt and sugar at a fixed value. Finally, constraints (\ref{eq:model_food_cons6}) and (\ref{eq:model_food_cons7}) characterize the CCL framework. In this sense, constraint (\ref{eq:model_food_cons6}) represents the learned decision variables: the food basket palatability, $y_{\text{pltb}}$. These score values are completely dependent from the basket components $\bm{x}$, establishing a relationship learned by means a predictive model $f^{\mathcal{D}}$, trained using a dataset $\mathcal{D}$. On the other hand, constraint (\ref{eq:model_food_cons7}) represents a chance constraint: we force the probability of the palatability score $y_{\text{pltb}}$ to be greater than $t$ at confidence level $\alpha$.

\subsubsection{Methodology}

As has been mentioned before, a dataset is needed to relate the palatability score with the different commodities. In this set, a sample dataset of 5000 observations has been employed, where each food basket is rated with a palatability value from zero (worst case) to one (best case). In this dataset, we can find up to 25 features or commodities, including beans, lentils, sugar, oil, and wheat, among others. 

For the sake of brevity, we refer the reader to \cite{maragno2021mixed}, where the authors explain the complete development of the dataset. In the same work, information about nutritional commodity contents and nutrient requirements can be found.

In relation to the predictive model, and again for illustrative purposes, a two hidden-layer deep neural network with 100 neurons per layer has been selected as the main function to establish the relationship between palatability and commodities. The final layer of this NN will be composed of only one neuron in the point and quantile methodology. In order to approximate the superquantile, we output 5 different quantiles in the final layer (see Section \ref{sec:superq_estim} for more information). Furthermore, a palatability lower bound of $0.5$ is chosen.

As in the previous case study, three different optimal solutions will be obtained for each point, quantile, and superquantile estimation methodology. For the point estimation methodology, we will assume that the chance constraint does not exist (that is, the optimization problem will only include a point prediction learned constraint). For the quantile and superquantile methodology, we force the chance constraint to be fulfilled at a confidence level of 95\%.

\subsubsection{Results}

The same software presented in Section \ref{sec:concrete_results} was employed to solve the palatable food basket problem. Table \ref{tab:food_sol} shows the optimal solution for the palatable food basket regarding the employed methodology.

\begin{table}[!ht]
    \centering
    \caption{Optimal food basket}
    \label{tab:food_sol}
    \begin{tabular}{c|ccc}
        \multicolumn{1}{p{4cm}|}{\centering \textbf{Decision Variable \\ Selection} } & \multicolumn{1}{p{3cm}}{\centering \textbf{Point} \\ \textbf{Estimation} } & \multicolumn{1}{p{3cm}}{\centering \textbf{Quantile} \\ \textbf{Estimation} } & \multicolumn{1}{p{3cm}}{\centering \textbf{Superquantile} \\ \textbf{Estimation} } \\ \hline
        Milk & 45.74 g & 45.55 g & 44.85 g  \\ 
        Salt & 5.00 g & 5.00 g & 5.00 g  \\
        Lentils & 32.87 g & 34.16 g & 38.97 g  \\
        Maize & 99.94 g & 98.53 g & 93.28 g \\ 
        Sugar & 20.00 g & 20.00 g & 20.00 g \\ 
        Oil & 21.66 g & 21.70 g & 21.86 g \\ 
        Wheat & 280.54 g & 280.68 g & 281.19 g  \\ 
        WSB & 74.69 g & 74.79 g & 75.18 g \\ \hline
        \textbf{Cost} & 3258.90 \$ & 3260.78 \$ & 3267.78 \$ \\ \hline
    \end{tabular}
\end{table}

For clarity, Table \ref{tab:food_sol} only contains the optimal solution for all the components of the food basket with a value greater than zero. As can be seen, differences across methodologies are minimum. We can notice how, as we move from point to quantile, and from quantile to superquantile estimation, the quantity of milk in the food basket decreases. The same happens with the amount of maize. On the other hand, the amount of lentils, oil, wheat, and WSB increases.

However, these small changes in the optimal food basket have an influence on the objective function. For the point estimation methodology, the total cost for fulfilling the demand takes a value of 3258.90\$. This cost rises up to 3260.78\$ with the quantile estimation methodology and reaches 3267.78\$ with the superquantile methodology.

The low variability in the three different food baskets also gives us information about how powerful Deep Neural Networks (DNN) methods are in predicting palatability. This accurate performance allows big chances in the probability distribution with small changes in the optimal solution. As a piece of additional information, the DNN obtained a testing MSE of $0.0001$ for the point estimation methodology, and a testing $\alpha$-quantile loss of $0.0013$ for the quantile methodology, confirming the excellent performance of the method.

\begin{figure}[!ht]
    \centering
    \includegraphics[width=\textwidth]{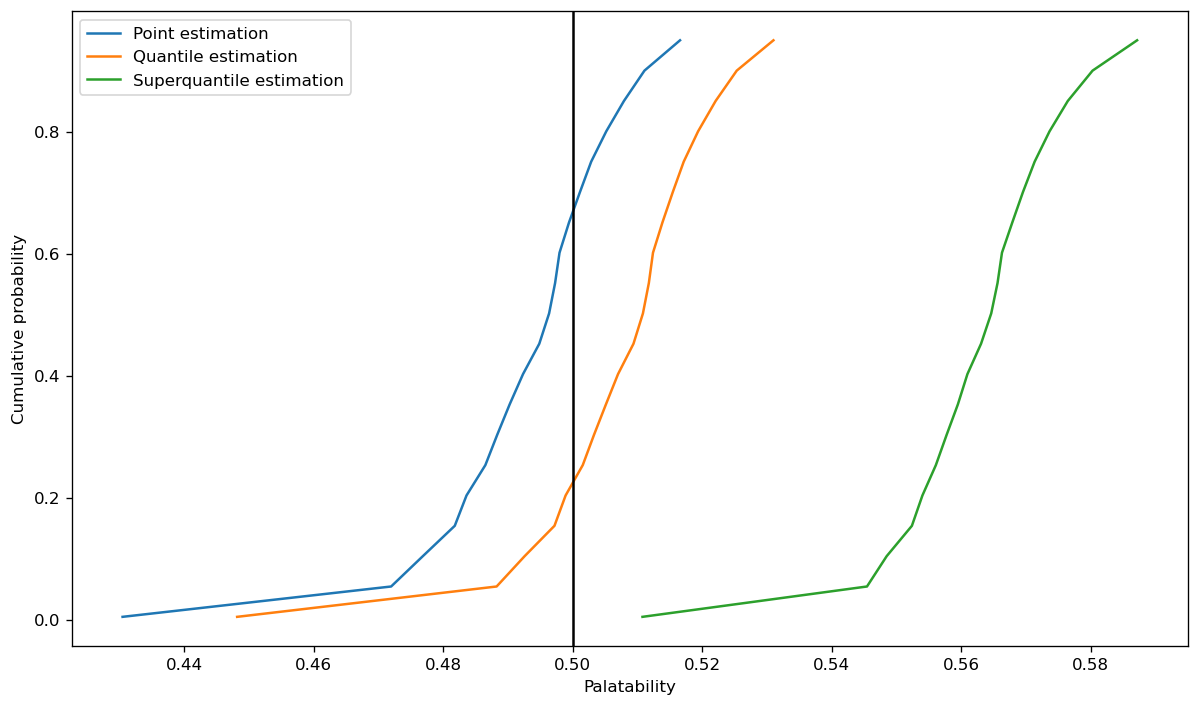}
    \caption{Optimal solution's distribution function for the palatable food basket problem regarding the employed methodology.}
    \label{fig:food_sol_dist}
\end{figure}

As in the previous case study, we present the distribution function of the optimal solution in Figure \ref{fig:food_sol_dist}. This distribution function has been computed by training another two hidden-layer NN, this time obtaining as outputs several conditional quantiles. Optimal solutions are passed through the fitted DNN and estimations of the quantiles are obtained to represent the distribution function.

We can see how, as we expected, the employment of chance constraints makes the distribution function of the optimal solutions to be over the seek palatability score ($0.5$) in most of the cases. We can also notice how powerful the DNN is by looking at the range of values of the palatability score. For example, the distribution function for the superquantile optimal solution ranges from $0.51$ to almost $0.59$. This narrow range actually represents a $99\%$ prediction interval for the estimated palatability. Note that this improvement is achieved by a basket cost increase of only $0.27\%$, which illustrates the convenience of our approach.

Finally, Table \ref{tab:food_perf} shows a brief summary of the performance of some selected predictive models in relation to the presented case study.

\begin{table}[!ht]
    \centering
    \caption{Chance constraint learning performance regarding methodology and predictive model for the palatable food basket problem}
    \label{tab:food_perf}
    \resizebox{\textwidth}{!}{%
    \begin{tabular}{l|c|cccccc}
        \multicolumn{1}{p{3cm}|}{\centering \textbf{Methodology} } & \multicolumn{1}{p{2cm}|}{\centering \textbf{Predictive model} } & \multicolumn{1}{p{4cm}}{\centering \textbf{Testing error (MAE/Quant loss)} } & \multicolumn{1}{p{2.5cm}}{\centering \textbf{Fitting \& CL time (s)} } & \multicolumn{1}{p{3.5cm}}{\centering \textbf{Optimization solving time (s)} } & \multicolumn{1}{p{3cm}}{\centering \textbf{Problem size (constraints)} } & \multicolumn{1}{p{2cm}}{\centering \textbf{Binary vars} } & \multicolumn{1}{p{2.25cm}}{\centering \textbf{Non-binary vars}} \\ \hline
        \multirowcell{3}{Point\\Estimation}  & LR & 0.120 & 0.2 & 0.2 & 5 & 0 & 26\\
        ~ & RF & 0.077 & 119.7 & 14.1 & 115 & 5574 & 76\\ 
        ~ & NN & 0.009 & 427.9 & 198.9 & 605 & 200 & 226\\ \hline
        \multirowcell{3}{Quantile\\Estimation} & LR & 0.014 & 1.3 & 0.3 & 5 & 0 & 26\\ 
        ~ & RF & 0.010 & 749.6 & 722.3 & 455 & 16651 & 176\\
        ~ & NN & 0.002 & 451.6 & 4.1 & 605 & 200 & 226\\ \hline
        \multirowcell{3}{Superquantile\\Estimation} & LR & 0.015 & 1.4 & 0.2 & 5 & 0 & 26\\
        ~ & RF & 0.010 & 739.9 & 477.1 & 455 & 16651 & 176\\ 
        ~ & NN & 0.002 & 515.6 & 10.2 & 610 & 200 & 231\\ \hline
    \end{tabular}%
    }
\end{table}

As can be seen in Table \ref{tab:food_perf}, the optimization solving time is directly (but not only) related to the problem size. We can take as an example the results produced from the different NN models. When employing a point estimation methodology, optimization time is way bigger than using (super)quantile estimation methodology. This may indicate that, for this particular problem, modelling the expected palatability value is more difficult than modelling a conditional quantile. However, when modelling palatability with RF models, it is faster to obtain the optimal food basket for point estimation methodology compared to the quantile ones, mainly due to the significant increment in the problem size (number of binary variables).

The accurate performance of NN-based models (especially with large datasets), jointly with their optimization solving time makes them suitable for these types of applications. As was shown before, NNs were accurate in correctly fitting within the CCL methodology.

\section{Conclusions}
\label{sec:conclu}

In this work, a novel methodology has been developed in order to deal with the intersection of chance constraint and constraint learning within mixed-integer optimization problems, two topics that are increasingly getting attention among the operations research community.

From the constraint learning side, we deal with problems where some optimization variables $x$ are directly controlled by the decision-makers, whereas other variables $y$ can be indirectly influenced as they depend on the decisions over $x$ and external covariates known when decisions are made. Recent approaches have been developed in the literature, allowing practitioners to embed linearizable machine learning methods within their optimization problem to estimate the value of $y$ conditioned to the values of $x$. However, there are many real-world applications where we need to set a lower (or upper) bound to this learned variable, and the fulfillment of this constraint is not guaranteed, i.e., it depends on its probabilistic behaviour, which is further conditioned by the fitting process of the machine learning model.

For that reason, we take ideas from chance constraints and further impose probabilistic guarantees, i.e., we want the learned variable to be over (or under) a pre-defined bound with a certain probability.  This is done by employing quantile estimation linearizable models as the predictive functions to be embedded in the problem. These models have the capacity to output a conditional quantile at a specific level $\alpha$. Setting the value of this quantile over (or under) the bound will result in the fulfillment of the chance constraint. Furthermore, this procedure can be extended to estimate superquantiles, a better option if we seek coherent risk measures. This novel proposed framework has been denoted as Chance Constraint Learning (CCL) and can be employed for both solving chance constraints and adding statistical guarantees to learned constraints.

This machine learning-based approach will also allow us to avoid the employment of scenarios to model the chance constraint or the assumption of certain probability distributions. Thus, depending on the selected predictive model, we will reduce the number of additional binary variables that have to be employed.

Finally, the proposed CCL methodology has been tested in two different case studies, where we check the differences in the optimal solutions regarding the selected methodology and the impact of the conservatism of chance constraints.

We consider the proposed CCL framework to be useful in real-world problems where a decision variable has to be learned and the decision-maker is interested in giving robustness to its optimal solution. To this end, open-access software has been developed to give practitioners the appropriate tools to implement the CCL methodology within their optimization problems. Furthermore, the development of this framework paves the way for future research regarding quantile estimation and optimization problems. For example, we can tackle uncertainties by making use of prediction intervals composed of two learned quantiles.

\section*{Credit authorship contribution statement}

\textbf{Antonio Alc\'antara:} Conceptualization, Methodology, Validation, Investigation, Software, Writing - Original Draft. \textbf{Carlos Ruiz:} Conceptualization, Methodology, Validation, Investigation, Writing - Review \& Editing, Funding acquisition.

\section*{Declaration of competing interest}

The authors declare that they have no known competing financial interests or personal relationships that could have appeared to influence the work reported in this paper.

\section*{Acknowledgements}
The authors gratefully acknowledge the financial support from MCIN/AEI/10.13039/ 501100011033, project PID2020-116694GB-I00 and from the FPU grant (FPU20/00916).

%\printbibliography 

\bibliographystyle{elsarticle-num} 
\bibliography{biblio}

\end{document}